\documentclass{article}
\usepackage{biblatex}
\usepackage{graphicx} 
\usepackage{tikz-cd}
\usepackage{amsthm}
\usepackage{amsmath}
\usepackage{amsfonts}
\usepackage{indentfirst}
\usepackage{adjustbox}
\theoremstyle{definition}
\newtheorem{exmp}{Example}[section]

\newenvironment{proposition}[2][Proposition]{\begin{trivlist}
\item[\hskip \labelsep {\bfseries #1}\hskip \labelsep {\bfseries #2.}]}{\end{trivlist}}

\title{Non-trivial Odd Dimensional Equivariant Cohomology of Torus Orbits of Generic Points in $gr(2,k)$.}
\author{Timothy Paczynski }
\date{}
\begin{document}

\maketitle

\section{Introduction}
Given a grassmannian $gr(2,k)$ we can classify the various orbits of  points under the action of a torus which we can represent as a matrix with diagonal entries $t_1$ through $t_k$. This orbits of various points give rise to various toric varieties. In this work we will consider only those toric varieties which arise as orbits of generic points. For non-generic points, please see the follow up to this document by the same author.

The orbit of these generic points of gr(2,k) are in general not GKM spaces and it would be interesting to have a nice understanding of the equivariant cohomology rings under this particular torus action. To understand these equivariant cohomology rings we will first describe a procedure of blowing these generic toric varieties at their singular points so that the resulting space is GKM. Then we will use tools from commutative algebra to glean information about our original toric varieties from this blown up toric variety. 

In pursuit of this goal, it is first useful to lay down some notation. Let $X$ be the closure of the orbit of a generic point in the grassmannian $gr(2,k)$ under $T_k$ and the fixed points in the closure of the this action be called $A$. We can immediately compute that there are exactly as many fixed points as there are transpositions in $s_{k+1}$, and similarly we can enumerate the one dimensional orbits under this action counting that there are $\frac{k(k-1)}{2}$. Each of these fixed points in the closure of our orbit of our generic point is a singular point on the original variety. These fixed points comprise all of the singularities our variety as can easily be observed from the combinatorial facts about either the fan or the resulting moment polytope. 

It is often useful from the equivariant prospective to consider the images of the moment maps fo these varieties. In particular when our variety is GKM its moment graph allows one to algorithmically compute the equivariant cohomology of these varieties via work of Atiyah, Bott, and Goresky, Macpherson, and Kottiwitz. One useful way to do this explicitly is to use the plucker embedding of the toric variety into projective space and embed the torus we acted on the grassmannian with into the torus that naturally acts on the projective space we send it into. 

To resolve the singularities of the toric variety we can blow-up at these singular points. In the work that follows we will blow up at all of the singular points simultaneously. It will be useful to think of a decorated moment polytope (moment graph) of this blown up space. This resulting polytope of the blow up will become a copy of the orginal toric varieties moment polytope except that it will have its corners cut off by hyperplanes which contain a point in the interior of moment polytope that is within an $\epsilon$'s distance form our original singular point and whose normal vector is just the vector defined by the difference of the singular point in question and the central point of our moment polytope. 

With a rough sketch of the objects at play in hand we can now begin the formal work.

\section{Preliminaries}

This paper will make extensive use about various facts widely used in the study of equivariant cohomology, toric variates, and homological/commutative algebra. For basic facts about toric varieties \cite{Fulton}, and for more detailed facts and some concrete computations in equivariant cohomology of the singular case \cite{CLS}. For a wonderful and down to earth introduction to equivariant cohomology a great reference is \cite{Tymo}, and a nice introduction which is grounded in geometry is \cite{Tu}. For more generic equivariant cohomology facts, a wonderful reference is \cite{GKM}. 

Let $X$ be the closure of the orbit of a generic point in the grassmannian $gr(2,k)$ under the action of $T_k$ and let A be the set of fixed points in the closure of this action. Since equivariant cohomology is a general cohomology theory we have a well defined cohomology of the pair $(X,A)$ which comes from the injection of $A$ into $X$. Thus we have long exact sequences of the form $$\cdots \rightarrow H^i_{T_k}(X) \rightarrow H^i_{T_k}(A) \rightarrow H^{i+1}_{T_k}(X,A) \rightarrow \cdots.$$

Let $\widetilde{X}$ denote the blow up of $X$ at the set of fixed points $A$ and let $\widetilde{A}$ denote the exceptional divisor of this blowup. Since there is a blow down map from $\widetilde{X}$ to $X$ and equivariant cohomology is a contravariant functor we have an induced map from the long exact sequences of the pair $(X,A)$ to the long exact sequence of the pair $(\widetilde{X}, \widetilde{A})$ as follows:

\begin{tikzcd}[cells={nodes={minimum height=2em}}]
\cdots \arrow[r] &  H^i_{T_k}(X) \arrow[r,"\phi"] \arrow[d, "f_X"] & H^i_{T_k}(A) \arrow[d,"f_A"] \arrow[r, "\delta"] &  H^{i+1}_{T_k}(X,A) \arrow[r]\arrow[d,"f_{(X,A)}"] \arrow[r, "\rho"] & \cdots\\
\cdots  \arrow[r] &  H^i_{T_k}(\widetilde{X}) \arrow[r,"\widetilde{\phi}"] &  H^i_{T_k}(\widetilde{A}) \arrow[r,"\widetilde{\delta}"] &  H^{i+1}_{T_k}(\widetilde{X},\widetilde{A}) \arrow[r, "\widetilde{\rho}"] & \cdots
\end{tikzcd}

It can immediately be noted that the space $(X,A)$ is equivalent to the space $(\widetilde{X}, \widetilde{A})$ since the only difference in between these spaces comes from blowing down the space $\widetilde{A}$ to points of $A$.

\begin{tikzcd}[cells={nodes={minimum height=2em}}]
\cdots \arrow[r] &  H^i_{T_k}(X) \arrow[r,"\phi"] \arrow[d, "f_X"] & H^i_{T_k}(A) \arrow[d,"f_A"] \arrow[r, "\delta"] &  H^{i+1}_{T_k}(X,A) \arrow[r]\arrow[d, "Id"] \arrow[r, "\rho"] & \cdots\\
\cdots  \arrow[r] &  H^i_{T_k}(\widetilde{X}) \arrow[r,"\widetilde{\phi}"] &  H^i_{T_k}(\widetilde{A}) \arrow[r,"\widetilde{\delta}"] &  H^{i+1}_{T_k}(\widetilde{X},\widetilde{A}) \arrow[r, "\widetilde{\rho}"] & \cdots
\end{tikzcd}

Next to examine these two sequences further we first need to prove that $\widetilde{X}$ and $\widetilde{A}$ are GKM spaces under the action of $T_k$. This follows from the fact that toric varieties admit resolutions of singularities to simplicial spaces and that simplicial spaces are GKM for equivariant cohomology with $Q$ coefficients. \cite{CLS} \footnote{ $\widetilde{X}$ is a rationally smooth toric variety and $\widetilde{A}$ is isomorphic to the disjoint union of $k(k-1)/2$ copies of $P^{k-3}\times P_1$ or a line bundle over $P^{k-3}$. }
Thus we know that the odd  dimensional equivariant cohomology of these spaces under this action vanishes when considering $\mathbb{Q}$ coefficients. Similarly, it is well know that the fixed point set of a torus action has no odd equivariant cohomology. Given these two facts our updated exact sequence appears as:

\begin{center}
\adjustbox{scale=.9,center}{
\begin{tikzcd}[cells={nodes={minimum height=2em}}]
 H^{2i}_{T_k}(X,A) \arrow[r]\arrow[d, "Id"] \arrow[r, "\rho"] & H^{2i}_{T_k}(X) \arrow[r,"\phi"] \arrow[d, "f_X"] & H^{2i}_{T_k}(A) \arrow[d,"f_A"] \arrow[r, "\delta"] &  H^{2i+1}_{T_k}(X,A) \arrow[r]\arrow[d, "Id"] \arrow[r, "\rho"] 
&  H^{2i+1}_{T_k}(X) \arrow[r,"\phi"] \arrow[d, "f_X"] & 0 \arrow[d,"f_A"] \\
 H^{2i}_{T_k}(\widetilde{X},\widetilde{A}) \arrow[r, "\widetilde{\rho}"] & H^{2i}_{T_k}(\widetilde{X}) \arrow[r,"\widetilde{\phi}"] &  H^{2i}_{T_k}(\widetilde{A}) \arrow[r,"\widetilde{\delta}"] &  H^{2i+1}_{T_k}(\widetilde{X},\widetilde{A}) \arrow[r, "\widetilde{\rho}"] & 0 \arrow[r,"\widetilde{\phi}"] &  0 
\end{tikzcd}
}
\end{center}

Next I would like to show that the map $\widetilde{\phi}$ is injective. Since both $\widetilde{A}$ and $\widetilde{X}$ are GKM, and the moment graph of $\widetilde{A}$ is a subgraph of the moment graph of $\widetilde{X}$ we get the injection $\widetilde{\phi}$ from the injection of the two spaces into their fixed point set. \footnote{It may be easier to observe that since $\widetilde{X}$ is GKM it injects into $\widetilde{A}$, and thus since $\widetilde{A}$ has no odd equivariant cohomology and since is a long exact sequence it must be the case that $(\widetilde{X},\widetilde{A})$ has no odd equivariant cohomology}

The injectivity of this map gives us that $ 0 =  H^{2k}_{T_k}(\widetilde{X}, \widetilde{A}) =  H^{2k}_{T_k}(X,A)$ since $H^{2k}_{T_k}(\widetilde{A})$ is $0$ because $\widetilde{A}$ is GKM. Thus we can once again update our exact sequence to get that.

\begin{center}
\begin{tikzcd}[cells={nodes={minimum height=2em}}]
 0 \arrow[hookrightarrow]{r}\arrow[d, "Id"] \arrow[r, "\rho"] & H^{2i}_{T_k}(X) \arrow[r,"\phi"] \arrow[d, "f_X"] & H^{2i}_{T_k}(A) \arrow[d,"f_A"] \arrow[r, "\delta"] &  H^{2i+1}_{T_k}(X,A) \arrow[r]\arrow[d, "Id"] \arrow[r, "\rho"] 
&  H^{2i+1}_{T_k}(X) \arrow[r,"\phi"] \arrow[d, "f_X"] & 0 \arrow[d,"f_A"] \\
 0  \arrow[hookrightarrow]{r}{\widetilde{\rho}} & H^{2i}_{T_k}(\widetilde{X}) \arrow[hookrightarrow]{r}{\widetilde{\phi}} &  H^{2i}_{T_k}(\widetilde{A}) \arrow[twoheadrightarrow]{r}{\widetilde{\delta}} &  H^{2i+1}_{T_k}(\widetilde{X},\widetilde{A}) \arrow[r, "\widetilde{\rho}"] & 0 \arrow[r,"\widetilde{\phi}"] &  0 
\end{tikzcd}    
\end{center}
and this further implies that our original $\phi$ was an injection. Thus even though our original space $X$ wasn't GKM it's equivariant cohomology modules in even dimensions still injection into the even equivariant cohomology modules of its fixed point set just as if it was GKM.

\begin{center}
\begin{tikzcd}[cells={nodes={minimum height=2em}}]
 0 \arrow[hookrightarrow]{r}\arrow[d, "Id"] \arrow[r, "\rho"] & H^{2i}_{T_k}(X) \arrow[hookrightarrow]{r}{\phi} \arrow[d, "f_X"] & H^{2i}_{T_k}(A) \arrow[d,"f_A"] \arrow[r, "\delta"] &  H^{2i+1}_{T_k}(X,A) \arrow[r]\arrow[d, "Id"] \arrow[r, "\rho"] 
&  H^{2i+1}_{T_k}(X) \arrow[r,"\phi"] \arrow[d, "f_X"] & 0 \arrow[d,"f_A"] \\
 0  \arrow[hookrightarrow]{r}{\widetilde{\rho}} & H^{2i}_{T_k}(\widetilde{X}) \arrow[hookrightarrow]{r}{\widetilde{\phi}} &  H^{2i}_{T_k}(\widetilde{A}) \arrow[twoheadrightarrow]{r}{\widetilde{\delta}} &  H^{2i+1}_{T_k}(\widetilde{X},\widetilde{A}) \arrow[r, "\widetilde{\rho}"] & 0 \arrow[r,"\widetilde{\phi}"] &  0 
\end{tikzcd}
\end{center}

This leads one to consider when $H^{2i+1}_{T_k}X = 0$, and more generally, if we can determine $H^{2i+1}_{T_k}(\widetilde{X}, \widetilde{A})$. 

In pursuit of the first question one notices, immediately from the fact that this is an exact sequences, that $$H^{2i+1}_{T_k}X \cong H^{2i+1}_{T_k}(X,A) \big/ \left( H^{2i}_{T_k}A / H^{2i}_{T_k} X) \right). $$
In general it is non-trivial to compute the quotients of free modules as above. Thus a first approach would be to try considering these spaces with $Q$ coefficients to get a dimension count on their differences.

\section{The Existence of Non-Trivial Odd Equivariant Cohomology}

Our first objective is to show when some of the odd dimensional equivariant cohomology of these modules is non-zero by doing a dimension count on the generators of these spaces in degree $2$ and showing that by difference of degree 2 generators it must be the case that $H^{3}_{T_k}(X) \neq 0$. 

To do this we need to find all the complex degree zero and degree 1 generators, which for $GKM$ spaces is equivalent to counting the number of vertices with $0$ and multiplying by $k$ and then adding the number of degree 1 vertices. 
First we can just observer that the number of vertices of $A$ is $k(k-1)/2$ and thus the number of degree two free generators of $H^2_{T_k}(A)$ is $k^2(k-1)/2$. In fact this implies that $Dim_\mathbb{Q}((H^{2n}_{T_k}(A)) =  \binom{n+k-1}{n}(k)(k-1)/2$.

Next, we can count the number of fixed points of $\widetilde{A}$. Since $\widetilde{A}$ is the disjoint union of $k(k-1)/2$ copies of moment graphs isomorphic to the moment graph of $P^{k-3} \times P^1$ we have $k(k-1)/2$ zero dimensional generators and $2k(k-1)/2$ degree 1 generators. Thus we have that the count on free generators of $H^2_{T_k}(\widetilde{A})$ is $k^2(k-1)/2 +k(k-1)$. In general $H^{2n}_{T_k}(\widetilde{A})$ is counted by: 

\begin{center}
\adjustbox{scale=1.15,center}{
$\binom{k+n-1}{n} \frac{(k)(k-1)}{2} + \sum^{k-3}_{i=1} \binom{k +n - i -1}{n - i} k(k-1)  + \binom{k+n -(k-2)-1}{n - (k-2)} \frac{(k)(k-1)}{2}.$}
\end{center}

This can be calculated by counting the possible new coefficients on the free generators, of which there are two in every dimension coming from every dimension except 0 and 2(k-2) which there is one. 

\begin{proposition}{1} Let $X$ be the closure of the orbit of a generic point in $gr(2,k)$ and let $\Tilde{X}$ be the blow up of $X$ which is a simplicial resolution with exceptional divisor $\tilde{A}$. Then:
    $$H^{2n}_{T_k}(\widetilde{A}: \mathbb{Q}) = \mathbb{Q}^{\binom{k+n-1}{n} \frac{(k)(k-1)}{2} + \sum^{k-3}_{i=1} \binom{k +n - i -1}{n - i} k(k-1)  + \binom{k+n -(k-2)-1}{n - (k-2)} \frac{(k)(k-1)}{2}}.$$
\end{proposition}

For the free generators of $H^2_{T_k}(\widetilde{X})$ we have to work slightly harder. 
To compute them we can notice that the moment graph has all the same vertexes as $\widetilde{A}$ above. It has the same edge set as well except for one additional directed edge either coming into or going out of each vertex. This is all determined form the blow up procedure and the direction comes from the original direction of the same edge in the moment graph of $\widetilde{X}$. Thus to compute the directed valence of any one vertex we just need to compute how this one new edge effects the two vertices it joins. 

We can visualize this by picking a vertex we want to consider the blow up of in the original polytope. If $(ij)$ is the vertex, then its blow up is the copy of $P^{k-3}$ defined by the vertices $(ik)$ where $k \neq j$ joined to the copy of $P^{k-3}$ defined by the vertices $(jk)$ via the edges between $(ik)$ and $(jk)$. This will be isomorphic to the moment graph of $P^{k-3} \times P^1$. We know the edges coming into each $(ik)$ is the original edge from the vertex $(ik)$ into $(ij)$. Thus if we represent the vertices of our blown up space as two arrays of length $n-2$ whose entries correspond to the flow along the vertices of $P^{n-3}$, we get that the first $j-2$ entries of $i$'s array increase by 1 and the first $i-1$ entries of $j's$ array increase by 1. The remaining entries stay the same. Since we know the original degree of the vertices in $i's$ array is $0,1,2, \cdots n-3$, and we know the degree of the vertices in $j's$ array is $1,2, \cdots, n-2$. Thus we can count the total vertices of degree $k$ by starting with $k(k-1)/2$
vertices of each of the types listed above (times two for $1$ through $n-3$) and then moving vertices from one bucket to the next one up when they have the additional edge raising their degree by 1. 

One way to see this is to notice that the edges going into the original point denoted $(ij)$ were all the edges from $(ik)$ with $k <j$ and all the edges form $(hj)$ with $h <i$. Since each of these edges now points to the vertex which represents itself via the slice we took, it immediately follows that we can count the change in these vertices degrees by simply recording the degrees of the vertices of $P^{k-1} \times P^1$, adding one to each which comes from a transposition lower than $ij$, and then summing over the total. 

An example of how to do this blow up for a point in $gr(2,3)$ is below where the green edges are the new edges we have added in and the red edge and its associated vertices are the blown up point(s). The original vertex $13$ has a blow up here with edge valence $(1,1) = (0+1, 1+0)$ 
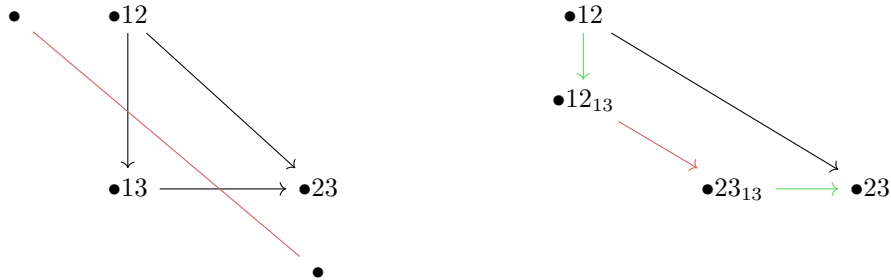
\begin{exmp}
\[\begin{tikzcd}
	\bullet & {\bullet 12} &&&&& {\bullet 12} \\
	&&&&&& {\bullet 12_{13}} \\
	& {\bullet 13} && {\bullet 23} &&&& {\bullet 23_{13}} & {\bullet 23} \\
	&&& \bullet
	\arrow[from=1-2, to=3-2]
	\arrow[from=1-2, to=3-4]
	\arrow[draw={rgb,255:red,92;green,214;blue,92}, from=1-7, to=2-7]
	\arrow[from=1-7, to=3-9]
	\arrow[draw={rgb,255:red,214;green,92;blue,92}, from=2-7, to=3-8]
	\arrow[from=3-2, to=3-4]
	\arrow[draw={rgb,255:red,92;green,214;blue,92}, from=3-8, to=3-9]
	\arrow[draw={rgb,255:red,214;green,92;blue,92}, no head, from=4-4, to=1-1]
\end{tikzcd}\]

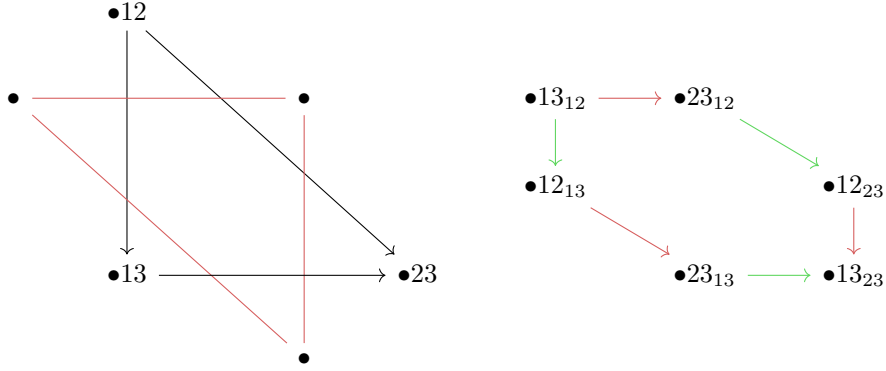
\[\begin{tikzcd}
	& {\bullet 12} \\
	\bullet &&& \bullet && {\bullet 13_{12}} & {\bullet 23_{12}} \\
	&&&&& {\bullet 12_{13}} && {\bullet 12_{23}} \\
	& {\bullet 13} &&& {\bullet 23} && {\bullet 23_{13}} & {\bullet 13_{23}} \\
	&&& \bullet
	\arrow[from=1-2, to=4-2]
	\arrow[from=1-2, to=4-5]
	\arrow[draw={rgb,255:red,214;green,92;blue,92}, no head, from=2-1, to=2-4]
	\arrow[draw={rgb,255:red,214;green,92;blue,92}, no head, from=2-1, to=5-4]
	\arrow[draw={rgb,255:red,214;green,92;blue,92}, no head, from=2-4, to=5-4]
	\arrow[draw={rgb,255:red,214;green,92;blue,92}, from=2-6, to=2-7]
	\arrow[draw={rgb,255:red,92;green,214;blue,92}, from=2-6, to=3-6]
	\arrow[draw={rgb,255:red,92;green,214;blue,92}, from=2-7, to=3-8]
	\arrow[color={rgb,255:red,214;green,92;blue,92}, from=3-6, to=4-7]
	\arrow[color={rgb,255:red,214;green,92;blue,92}, from=3-8, to=4-8]
	\arrow[from=4-2, to=4-5]
	\arrow[color={rgb,255:red,92;green,214;blue,92}, from=4-7, to=4-8]
\end{tikzcd}\]
\end{exmp}


To compute the degree of the vertices in the moment graph of $\widetilde{X}$ we need to consider how the valence of the vertices of the moment graph isomorphic to moment graph of $P^{k-3} \times P^1$ change when they gain an extra directed edge coming from the original moment graph of $X$. Thus let us adopt the convention that given $a$ and $b < c$ our edge is directed from $ab$ to $ac$. \footnote{This is possible in the first place since our the labels on our graph are a combinatorial morse function.} Then the graph $P_{ij}$, corresponding to the blown up point $ij$ before the additional edge added to each vertex coming from outside the polytope's valence is added, will have 1 valence 0 vertex, 2 valence 1 vertices and some other vertices valence 2 or more. 

Since no other vertices will generate any portion of $H^2(\widetilde{X})$ lets initially work with these. For a vertex of the form $(1j)$ with $j >3$ we have that $P_{1j}$ has degree 0 vertex $(12_{1j})$ and degree 1 vertices labeled $(13_{1j})$ and $(2j_{1j})$. Since $(12)$ and $(13)$ are below $1j$ we know these vertices valence will increase by 1. Thus, $(12)$ and $(2j)$ will be valence 1 in the moment graph of $\widetilde{X}$ while $(13)$ will be valence 2. When our vertex is of the form $(2j)$ with $j>3$, then our degree 0 vertex is $(12)$ and our degree 1 vertices are $(23)$ and $(1j)$. Since $(12), (1j)$, and $(23)$ are all less than $2j$, we have that $(12)$ will become valence 1 and $(23)$ and $(1j)$ will become valence 2. When our vertex is $(ij)$ with $i >2$, $i<j$, and $j>3$ then our degree 0 vertex is $(1i)$ and our degree 1 vertices are $(2i)$ and $(1j)$. So in this case $(1i)$ will become valence 1 and $(2i)$ and $(1j)$ will become valence 2. This captures all cases except for when the blown up vertex is $(12)$, $(13)$, and $(23)$. In these cases the vertices of valence 0 and 1 valence in $P_{ij}$ in the moment graph of $\widetilde{X}$ will become: $(12)$: 1 valence 0 vertex and 2 valence 1 vertices; $(13)$: 3 valence 1 vertices; $(23)$: 2 valence 1 vertices and 1 valence 2 vertex.

Thus we can compute the total number of valence 0 vertices in $\widetilde{X}$ is 1 and the number of degree 1 vertices in $\widetilde{X}$ for $gr(2,k)$ with $k>2$ is
 $$7 + 2(k-3)  +  ((k)(k-1)/2 -3).$$ Therefore the number of free generators of $H^2_{T_k}(\widetilde{X})$ when considered with $\mathbb{Q}$-coefficients is $$2k + (k)(k-1)/2 +1.$$
 Since $\widetilde{\phi}:H^2_{T_k}(\widetilde{X}) \rightarrow H^2_{T_k}(\widetilde{A})$ is injective this implies that 
 \begin{align*}
   rank(H^3_{T_k}(X,A))&= rank(H^3_{T_k}(\widetilde{X}, \widetilde{A}))   ) \\
   &=  k^2(k-1)/2 +k(k-1) - (2k + (k)(k-1)/2 +1)
 \end{align*}
for $gr(2,k)$ with $k >2.$ Using the computation of $H^2_{T_k}(A)$ above we have that the rank of $H^2_{T_k}(A)$ minus the rank of $H^3(X,A)$ is $$k(k-1) - (2k + (k)(k-1)/2 +1 = \frac{1}{2}(k^2 -5k-2).$$
This polynomial has roots $\frac{1}{2}(5 \pm \sqrt{33})$ and hence we have shown that $k>5$ the map $\delta: H^2_{T_k}(A) \rightarrow H^3_{T_k}(X,A)$ has cokernel of rank at least $\frac{1}{2}(k^2 -5k-2) >0$ and thus that there is odd equivariant cohomology for the  torus orbits of generic points of the grassmannians $gr(2,k)$ with $k>5$.
\newtheorem{theorem}{Theorem}
\begin{theorem}[] \label{thm:odd_cohomology}
The Equivariant Cohomology of the torus orbit of a generic point in the $gr(2,k)$ where $k >= 6$ has non-zero odd equivariant cohomology in degree $3$.
\end{theorem}
It is now natural to ask what happens in smaller grassmannians and if we can compute odd cohomology in these cases. For $Gr(2,5)$ we need to examine the group $H^4_{T_k}(\widetilde{X})$ and $H^4_{T_k}(\widetilde{A})$. In this case we (after ignoring the degree zero generators of $H^*_{T_k}\widetilde{A}$ which will cancel in the exact sequence with $H^*_{T_k}(A)$ that $H^4_{T_k}(\widetilde{A})$ has $ 120 = 5 *20 +20 $ free generators and $H^4_{T_k}(\widetilde{X})$ has $116 = 10+ 5 *16 + 26$ thus the odd dimensional cohomology shows up later when just examining rank in this case. For $Gr(2,4)$ rank won't be sufficient to show that odd equivariant cohomology exists so we will need something strong. However, as we will show later in this paper. The rank argument augmented by facts about the spectral sequence associated to the equivariant cohomology modules will let us determine that there is non-trivial odd cohomology. 
\section{Explicit Computations of $H^3_{T_k}$}

Recall, the Cartan Mixing Diagram used to define the spectral sequence associated for the equivariant cohomology of a space, $X$, with a $G$-action is:
\[\begin{tikzcd}
	EG && { EG \times X} && X \\
	\\
	BG && { EG \times_G X} && {X /G}.
	\arrow[from=1-1, to=3-1]
	\arrow[from=1-3, to=1-1]
	\arrow[from=1-3, to=1-5]
	\arrow[from=1-3, to=3-3]
	\arrow[from=1-5, to=3-5]
	\arrow[from=3-3, to=3-1]
	\arrow[from=3-3, to=3-5]
\end{tikzcd}\]

Since our $G$ here is a torus, $T_k$, we can rewrite our diagram with the specific base space and classifying spaces associated to $EG$ and $BG$,
\[\begin{tikzcd}
	\prod^k S^\infty && { (\prod^k S^\infty)\times X} && X \\
	\\
	\prod^k\mathbb{CP}^\infty && {(\prod^k S^\infty) \times_G X} && {X /G}.
	\arrow[from=1-1, to=3-1]
	\arrow[from=1-3, to=1-1]
	\arrow[from=1-3, to=1-5]
	\arrow[from=1-3, to=3-3]
	\arrow[from=1-5, to=3-5]
	\arrow[from=3-3, to=3-1]
	\arrow[from=3-3, to=3-5]
\end{tikzcd}\]
From now on we will also consider $X$ to be the orbit of a generic point by $T_k$ in $gr(2,k)$.

The key fact about these spectral sequence diagrams is that for reasonably nice spaces, like we have, there is a fibration formed with fiber $X$, total space $(\prod^k S^\infty) \times_G X$ and base $\prod^k\mathbb{CP}^\infty$. This fibration admits a spectral sequence for the total space whose $E^2$ page is the tensor product of copies of the cohomology modules of $X$ and $\prod^k\mathbb{CP}^\infty$ respectively.

This is a boon since the cohomology modules of $\prod^k\mathbb{CP}^\infty$ are well understood, and the cohomology modules of toric varieties are somewhat understood. In particular, we know the following the methods of \cite{Jor}, written succinctly in \cite{CLS}, one can show that a particular spectral sequence of filtered tori which compose a toric variety degenerates on the $E^2$-page if we use $\mathbb{Q}$-coefficients. While the details are fairly compelling, the main point is that the values of $E^{p,q}_1$ are computable via combinatorial methods, and that it is first quadrant spectral  sequence. 

To compute the $E^2$-page of the spectral sequence once has to count cones of the fan of the toric variety one is interested in. Explicitly,\footnote{this isn't actually $E^1$, it is a complex used to compute $E^2$ of the spectral sequence we care about that comes from a filtration via tori. I should  make that more evident.} $$E_1^{p,q}= \bigoplus_{\tau \in \Sigma(n-p)} \wedge^q M(\tau)  \cong  Q^{{p \choose q} |\Sigma(n-p)|}$$
where $\Sigma(k)$ is the number of cones of dimension $k$, $M$ is the dual lattice of the toric variety, and $M(\tau) := \tau^\perp \cap M$.

To define the differentials we first need to assign an orientation to every cone in $\Sigma$. Once we have done this the differentials on the $E^1$-page are, $$\delta^p:\bigoplus_{\tau \in \Sigma(n-p)} \wedge^q M(\tau) \rightarrow \bigoplus_{\tau \in \Sigma(n-(p+1))} \wedge^q M(\tau),$$
which are defined by component wise on the cones $(\tau, \sigma)$ in the direct sums by the rule that: if $\sigma \not\prec\tau$, then the component is $0$. Otherwise take any vector $v \in \tau, v \not\in \sigma$. Then if the vectors of  $\sigma$ and $v$ induce the same orientation as we originally selected for $\tau$ we say the orientation coefficient $c$ on the on the induced map from $i_{\tau,\sigma}^q:\wedge M(\tau) \rightarrow M(\sigma)$ is 1. If the orientation is reversed we say $c$ is negative 1.

To compute the differentials on the $E^1$-page with integral coefficients can be fairly tedious since one needs to consider the kernels of integer maps. Even though our vectors for these toric varieties fans are relatively simple, this can be computationally complex. Instead, we can side step this difficulty by first noting that for our complete toric varieties, on the $E_2$ page we will get a copy of $\mathbb{Z}$ and then all zeros for the index 0 row.  Thus we can just compute $H^3_T$ by just determining the maps between $E^{1,1}_1$ and $E^{2,1}_1$ and $E^{3,1}$. In particular this implies that $E^{2,2}_2$ is $H^2_T(X)$ and thus adding $\frac{1}{2}(k^2 -5k -2)$ determines $H^3_T(X)$. 

To see an explicit example, the case where our toric variety arises as the orbit of a generic point in $gr(2,k)$ with $k >= 4$ we can identify that these modules are $\wedge^1 M(\Sigma(1))$, $\wedge^1 M(\Sigma(2))$ and $\wedge^1 M(\Sigma(3))$. The maps between them are defined by how these cones map forward. The resulting cohomology of this toric variety with  $\mathbb{Q}$-coefficients is $(Q^1,0,Q^1,Q^2,Q^5,0,Q^1)$. We can now make use of spectral sequence associated to the equivariant cohomology of this orbit where along one side we have the toric variety and along the other side we have the torus we are working with $T^4$. A truncated version of the $E^2$ page looks as follows:

\[\begin{tikzcd}
	{\mathbb{Q}^5} &&&& \\
	{\mathbb{Q}^2} & 0 \\
	{\mathbb{Q}} & 0 & {\mathbb{Q}^4} \\
	0 & 0 & 0 & 0 & {} \\
	{\mathbb{Q}} & 0 & {\mathbb{Q}^4} & 0 & {\mathbb{Q}^{10}}
\end{tikzcd}\]

We can now clearly see that the following cohomology modules are already determined, $H_{T_4}^0(X) =\mathbb{Q}$, $H_{T_4}^1(X) = 0$, and $H_{T_4}^2(X) =\mathbb{Q}^5$. We can now simply compute that the map from $H_{T_4}^2(X) \rightarrow H_{T_4}^2(A)$ leaves us with cokernel of dimension $5$. we can compute the dimension of the cokernel of this map into $H_{T_k}(X,A)$ for the first section and get that its cokernel should be of dimension $5 - 3 = 2$. Thus from the $E^2$-page we can immediately deduce that $H_{T_4}^3 = \mathbb{Q}^2$.

Note that there is nothing stopping one from producing similar computations for the other torus orbits in $gr(2,k)$ except that the author has not worked out a nice closed form formula to compute the relevant $E^{i,j}$'s of the relevant spectral sequences for $i+j \in [3]$. Apriori, the same techniques hold for all computations of $H_{T_k}^3(X)$, though one could potentially get a lower bound which was not tight. 

\section{Conclusion}
In this paper, a family of torus orbit closures of generic points in $gr(2,k)$ have been introduced. It was supposed that these orbits should have non-trivial odd equivariant cohomology and we have produced a witness for such behavior. We have also introduce techniques for computing (or both upper and lower bounding) the cohomology modules in dimension 3. Similar techniques for lower bounding will hold for higher dimensions by simply extending the counting arguments from generators of dimension 0 and 1 above to generators of the blowup of higher dimensions as well.

\clearpage
\printbibliography 

\end{document}